\documentclass[10pt,twoside]{article}

\usepackage{amsmath}
\usepackage{amsfonts}
\usepackage{amssymb}
\usepackage{amscd}
\usepackage{amsthm}
\usepackage{amsbsy}
\usepackage{graphicx}
\usepackage{bm}

\def\@oddhead{\hfill \shorttitle \hfill \thepage}
\def\@evenhead{\thepage \hfill \shortauthor \hfill}
\def\@oddfoot{}
\def\@evenfoot{}

\newtheorem{thm}{Theorem}[section]
\newtheorem{defn}{Definition}[section]
\newtheorem{prob}[thm]{Problem}
\newtheorem{conj}[thm]{Conjecture}

\date{}
\title{\ \\[0.4cm] \ \\ \bf  Chern conjecture and isoparametric hypersurfaces\footnote{The project is partially supported
by the NSFC ( No.11071018 and No.11001016 ), the SRFDP, and the
Program for Changjiang Scholars and Innovative Research Team in
University.} }
\author{Jianquan Ge\footnote{School of Mathematical Sciences, Laboratory of Mathematics and Complex Systems, Beijing Normal
University, Beijing 100875. E-mail: jqge@bnu.edu.cn}\hspace{2mm} and
Zizhou Tang\footnote{The corresponding author. School of
Mathematical Sciences, Laboratory of Mathematics and Complex
Systems, Beijing Normal University, Beijing 100875. E-mail:
zztang@mx.cei.gov.cn}}

\begin{document}

\maketitle


\thispagestyle{empty}

\begin{abstract}
\vskip 3mm\footnotesize{

\vskip 4.5mm \noindent In this note we will review the most
important results and questions related to Chern conjecture and
isoparametric hypersurfaces, as well as their interactions and
applications to various aspects in mathematics.

\vspace*{2mm} \noindent{\bf 2000 Mathematics Subject Classification:
53-02, 53C40.}

\vspace*{2mm} \noindent{\bf Keywords and Phrases: minimal
submanifold, isoparametric hypersurface, Simons inequality, DDVV
conjecture.}}

\end{abstract}

\section{Introduction}
In this note we will review the most important results and questions
related to Chern conjecture and isoparametric hypersurfaces, as well
as their interactions and applications to various aspects in
mathematics.

We will start with a brief history of Chern conjecture and its
generalizations in Section \ref{chern conj}, then we introduce
isoparametric hypersurfaces in real space forms, complex space
forms, projective spaces, symmetric spaces, general Riemannian
manifolds and in particular, exotic spheres in Section \ref{isop},
and finally in Section \ref{int-app}, we discuss the interactions
and applications of Chern conjecture and isoparametric
hypersurfaces.

The exposition here is not rather complete, since our emphasis is on
open problems and possible future research interests instead of
historical results. For a more detailed history of results on the
``hypersurface" part of Chern conjecture and for an excellent survey
on isoparametric hypersurfaces and their generalizations, one can
see \cite{SW08} and \cite{Th00} (or \cite{Ce} for an updated one),
respectively.


\section{Chern conjecture}\label{chern conj}
More than 40 years ago, S.S. Chern, the leader in modern
differential geometry, proposed the following problem in several
places (cf. \cite{chern},\cite{C71},\cite{CCK}).
\begin{prob}\label{chern prob}
Consider closed minimal submanifolds $M^n$ immersed in the unit
sphere $S^{n+m}$ with second fundamental form of constant length
whose square is denoted by $\sigma$. Is the set of values for
$\sigma$ discrete? What is the infimum of these values of
$\sigma>\frac{n}{2-\frac{1}{m}}$?
\end{prob}
The affirmative hand of the first question, \emph{i.e.}, the set of
values for $\sigma$ should be discrete, is usually called Chern
conjecture. Up to now it is still far from a complete solution of
this problem even in the case when $M$ is a hypersurface (cf.
Problem 105 in \cite{Y82}). Moreover, as the advances of studies of
isoparametric hypersurfaces, people turned the hypersurface case
into the following new formulation (cf. \cite{V86}, \cite{SW08}).
\begin{conj}\label{ch-isop-conj}
Let $M^n$ be a closed, minimally immersed hypersurface of the unit
sphere $S^{n+1}$ with constant scalar curvature. Then $M$ is
isoparametric.
\end{conj}
This formulation bases on the fact that all known examples of
minimal hypersurfaces with constant scalar curvature (which, by
minimality, is equivalent to the condition that the second
fundamental form has constant length) in $S^{n+1}$ are
isoparametric, \emph{i.e.}, all of their principal curvatures are
constant, and $\sigma$ for these isoparametric hypersurfaces could
only attain finite number of values (see Sections below).

Chern conjecture originated from the famous Simons pinching theorem
(\cite{Simons}):
\begin{thm}
Let $M^n$ be a closed, minimally immersed submanifold in the unit
sphere $S^{n+m}$ and $\sigma$ the squared norm of its second
fundamental form. Then
$$\int_M \sigma\left(\sigma-\frac{n}{2-\frac{1}{m}}\right) dM\geq0.$$
 In particular, for $\sigma\leq\frac{n}{2-\frac{1}{m}}$
one has either $\sigma=0$ or  $\sigma=\frac{n}{2-\frac{1}{m}}$
identically on $M$.
\end{thm}
This theorem tells that if $\sigma$ is constant, it can not take any
value in the open interval $(0, \frac{n}{2-\frac{1}{m}})$. By his
always elegant presentation of calculations with moving frame
method, Chern \cite{chern}
 recovered the inequality of Simons above and further he (joint with do Carmo and Kobayashi,
 cf. \cite{CCK}, see also \cite{L}) classified the submanifolds on which
 $\sigma\equiv\frac{n}{2-\frac{1}{m}}$ to be only the Clifford
 minimal hypersurfaces and the Veronese surface (in $S^4$). Based on
 these, he proposed Problem \ref{chern prob}. After that, to this
 problem of submanifold part, \emph{i.e.}, for submanifolds of high codimension, many contributions had been made and could be
 concluded in the following theorem (cf. \cite{Y75}, \cite{S89}, \cite{CX93},
 \emph{etc.}).
 \begin{thm}\label{pinching-sub}
Let $M^n$ be a closed, minimally immersed submanifold in the unit
sphere $S^{n+m}$, $m\geq2$, $\sigma\leq\frac{2n}{3}$ everywhere on
$M$. Then $M$ is either a totally geodesic submanifold or a Veronese
surface in $S^{2+m}$.
 \end{thm}
On the other hand, to the hypersurface part of Chern conjecture,
Peng and Terng \cite{Peng-T1} made the first remarkable
breakthrough\footnote{Added in proof: Very recently Ding and Xin
\cite{DX} sharpened this result by neglecting the condition of
constant scalar curvature.}:
\begin{thm}\label{Peng-T}
Let $M^n$ be a closed, minimally immersed hypersurface in the unit
sphere $S^{n+1}$ with constant scalar curvature. Then there exists a
positive number $C(n)\geq\frac{1}{12n}$ such that if $\sigma>n$ then
$\sigma>n+C(n)$. Furthermore, if $n=3$ and $\sigma>3$, then
$\sigma\geq 6$.
\end{thm}
From then on, except for $3$-dimensional hypersurfaces there's no
more essentially affirmative answer to Chern conjecture though
considerable improvements of Theorem \ref{Peng-T} were made by
several geometers. For instance, Yang and Cheng (See, for example,
\cite{YC}) sharpened the constant $C(n)$ to $C(n)\geq
\frac{26}{61}n-\frac{16}{61}>\frac{1}{3}n$ and, under the additional
assumption that the sum of cubes of the principal curvatures $f_3$
is constant, $C(n)\geq\frac{13}{15}n-\frac{4}{5}\geq\frac{2}{3}n$.
Following \cite{Peng-T2}, Chang \cite{C3} showed that if the
hypersurface has exactly three pairwise distinct principal
curvatures in every point, then it is isoparametric. As for the
$3$-dimensional hypersurface case, combining results of
Almeida-Brito \cite{AB2} and Chang \cite{C1}, Conjecture
\ref{ch-isop-conj} was proved to be right even in a more general
category:
\begin{thm}\label{chern3dim}
Let $M^3$ be a closed hypersurface immersed in the unit sphere $S^4$
with constant mean curvature and constant scalar curvature. Then $M$
is isoparametric.
\end{thm}
Furthermore, the classification of $3$-dimensional hypersurfaces in
$S^4$ with two constant mean curvature functions \footnote{There're
totally $3$ mean curvature functions as elementary symmetric
polynomials of principal curvatures.} was also established and only
one case (vanishing mean curvature and Gauss-Kronecker curvature)
there admits non-isoparametric examples (cf.
\cite{AB1},\cite{AB2},\cite{AB3},\cite{C1},\cite{LO},\cite{R}). See
\cite{SW08} for an conclusion of these results, where it mentioned
that Bryant conjectured there should be a local version of Theorem
\ref{chern3dim}, \emph{i.e.}, minimal hypersurfaces in $S^4$ with
constant scalar curvature are isoparametric, and Chang \cite{C2}
proved some partial results towards this problem.

As the study of relations of intrinsic invariants and extrinsic
invariants goes deeper, it turns out that there're many more
Simons-type inequalities for submanifolds in space forms even in the
local sense which would certainly get powerful in attacking global
problems like Chern conjecture. In 1999, De Smet, Dillen,
Verstraelen and Vrancken \cite{DDVV} obtained an inequality
involving the (normalized) scalar curvature $\rho$, normal scalar
curvature $\rho^{\bot}$, and norm of mean curvature vector field $H$
for submanifolds of codimension $2$ in space forms $N(c)$ of
constant sectional curvature $c$, namely
 \begin{equation}\label{DDVV}
 \rho+\rho^{\bot}\leq|H|^2+c.
 \end{equation}
Then they proposed the conjecture that the inequality (\ref{DDVV})
holds for submanifolds of arbitrary codimension in space forms (the
so-called DDVV conjecture, and the inequality (\ref{DDVV}) is now
called the DDVV inequality). Recently, the DDVV conjecture was
completely solved by Lu \cite{L08} and the authors \cite{GT1} (where
we obtained further the pointwise equality condition) independently
with rather different methods (See \cite{GT2} for a survey). In
particular, by algebraic inequalities consisting commutators of
shape operators (in matrix form) that derived in the proof of the
DDVV inequality (\ref{DDVV}), Lu \cite{L08} generalized the results
of \cite{CCK} and Theorem \ref{pinching-sub} to the following
pinching theorem.
\begin{thm}\label{lu pinching}
Let $M^n$ be a closed, minimally immersed submanifold in the unit
sphere $S^{n+m}$. Let $\lambda_2$ be the second largest eigenvalue
of the semi-positive symmetric matrix $S:=(\langle
A^{\alpha},A^{\beta}\rangle)$ where $A^{\alpha}$s
($\alpha=1,\cdots,m$) are the shape operators of $M$ with respect to
a given (local) normal orthonormal frame. Suppose
$0\leq\sigma+\lambda_2\leq n.$
 Then $M$ is totally
geodesic, or is a Clifford minimal hypersurface in $S^{n+1}\subset
S^{n+m}$, or is a Veronese surface in $S^4\subset S^{2+m}$.
\end{thm}
Obviously, $\lambda_2\leq\frac{1}{2}\sigma$ and thus Theorem \ref{lu
pinching} is an elegant combination of those known pinching results
for hypersurfaces and submanifolds of high codimension. Based on
this, Lu \cite{L08} proposed a similar conjecture as Chern
conjecture by taking $\sigma+\lambda_2$ instead of $\sigma$:
\begin{conj}\label{lu conj}
Let $M^n$ be a closed, minimally immersed submanifold in the unit
sphere $S^{n+m}$ with constant $\sigma+\lambda_2$. If
$\sigma+\lambda_2>n$, then there is a constant $\varepsilon(n,m)>0$
such that $\sigma+\lambda_2>n+\varepsilon(n,m)$.
\end{conj}
Note that for hypersurface case it reduces to the corresponding part
of Chern conjecture and was proved by \cite{Peng-T1}. A difficulty
to attack this conjecture would be the non-smoothness of $\lambda_2$
in general.

\section{Isoparametric hypersurfaces}\label{isop}
The surveys on isoparametric hypersurfaces by Thorbergsson
\cite{Th00} and Cecil \cite{Ce} are so highly professional and
excellent that there's no reason for us to forget them and convince
ourself to reproduce a lengthy and boring story. Therefore, we would
like to pick out some important aspects to our own interests and
refer the details to \cite{Th00}, \cite{Ce}.

A (connected) hypersurface $M^n$ in a space form $N^{n+1}(c)$ of
constant sectional curvature $c$ is said to be \emph{isoparametric},
if it has constant principal curvatures. An isoparametric
hypersurface in $\mathbb{R}^{n+1}$ can have at most two distinct
principal curvatures, and it must be an open subset of a hyperplane,
hypersphere or a spherical cylinder $S^k\times \mathbb{R}^{n-k}$.
This was first proved for $n=2$ by Somigliana \cite{So} in 1919 (see
also Levi-Civita \cite{LeC}) and for arbitrary $n$ by Segre
\cite{Se} in 1938. A similar result holds in $\mathbb{H}^{n+1}$.
These two situations can be derived directly from the Cartan
identity which implies the number $g$ of distinct principal
curvatures must be $1$ or $2$ as \'{E}lie Cartan did in \cite{Ca1}
and \cite{Ca3}. However, as Cartan \cite{Ca1,Ca2,Ca3,Ca4} showed in
a series of four papers published in the period 1938--1940, the
theory of isoparametric hypersurfaces in the sphere $S^{n+1}$ is
much more beautiful and complicated. So far the classification
problem in this situation remains open for merely several cases,
although much progress has been made after the remarkable works of
M\"{u}nzner \cite{Mu80} which showed that the number $g$ of distinct
principal curvatures must be $1,2,3,4$, or $6$. In fact, the cases
when $g\leq3$ have been completely classified by Cartan. Besides,
Cartan and M\"{u}nzner found that each isoparametric hypersurface in
a sphere determines a so-called isoparametric function $f$
($=F|_{S^{n+1}}$, where $F: \mathbb{R}^{n+2}\rightarrow \mathbb{R}$
is called the \emph{Cartan polynomial}) that satisfies the so-called
Cartan-M\"{u}nzner equations:
\begin{eqnarray}
  &&|\nabla F|^2 =g^2|x|^{2g-2}, \label{eq1.1}\\
  &&\Delta F =\frac{g^2}{2}(m_2-m_1)|x|^{g-2}\label{eq1.2},
\end{eqnarray}
and conversely the set of level hypersurfaces of an isoparametric
function $f$ consists of a family of parallel hypersurfaces with
constant mean curvature (which in space forms implies each
hypersurface has constant principal curvatures) whose focal sets are
just the two critical sets $M_{\pm}$ of the isoparametric function
(as preimages of the maximum and minimum) and are submanifolds
(called the \emph{focal submanifolds}) of codimension $m_1+1$,
$m_2+1$ in $S^{n+1}$ respectively. In fact, M\"{u}nzner proved
further that each isoparametric hypersurface $M_t$ in the family
separates the sphere $S^{n+1}$ into two connected components $B_1$
and $B_2$, such that $B_1$ is a disk bundle with fibers of dimension
$m_1 + 1$ over $M_{+}$, and $B_2$ is a disk bundle with fibers of
dimension $m_2 + 1$ over $M_{-}$, \emph{i.e.}, they are tubes around
$M_{\pm}$. Moreover, the multiplicities of the distinct principal
curvatures $\kappa_1>\cdots>\kappa_g$ are alternatively $m_1$ and
$m_2$ ($m_1=m_2$ for $g=3$), and the principal curvatures of the
focal submanifolds $M_{\pm}$ in any normal direction are the
constants $\{\cot\frac{\pi}{g},\cdots,\cot\frac{(g-1)\pi}{g}\}$ with
multiplicities alternatively $m_2$ and $m_1$ for $M_{+}$
(alternatively $m_1$ and $m_2$ for $M_{-}$). Consequently, the focal
submanifolds are (austere) minimal submanifolds.

All known examples of isoparametric hypersurfaces with four
principal curvatures are of OT-FKM-type (examples with explicit
Cartan polynomials constructed by using representations of Clifford
algebras in \cite{OT76} and \cite{FKM}) with the exception of two
homogeneous families, having multiplicities $(m_1,m_2)$ (arranged to
be $m_1\leq m_2$) being (2, 2) and (4, 5). Beginning with
M\"{u}nzner , many mathematicians, including Abresch \cite{A}, Grove
and Halperin \cite{GH}, Tang \cite{Ta1} and Fang \cite{Fa}, found
restrictions on the multiplicities of the principal curvatures of an
isoparametric hypersurface with four or six principal curvatures.
This series of results culminated with the paper of Stolz \cite{St},
who proved that the multiplicities of an isoparametric hypersurface
with four principal curvatures must be the same as those in the
known examples of OT-FKM-type or the two homogeneous exceptions.
Cecil, Chi and Jensen \cite{CCJ} and independently Immervoll
\cite{Im} then showed that if the multiplicities $(m_1,m_2)$ of an
isoparametric hypersurface $M^n\subset S^{n+1}$ with four principal
curvatures satisfy $m_2\geq 2m_1-1$, then $M$ must be of
OT-FKM-type. Combining with results of Takagi \cite{Tak} and
Ozeki-Takeuchi \cite{OT76} that $M$ must be homogeneous if $m_1=1$
or $2$, their theorem classified isoparametric hypersurfaces with
four principal curvatures for all possible pairs of multiplicities
except for four cases, the homogeneous pair $(4,5)$, and the OT-FKM
pairs $(3,4)$, $(6,9)$ and $(7,8)$.  By employing more commutative
algebra than that explored in \cite{CCJ}, Chi \cite{Chi} gave a
proof to the case of $(3,4)$ that it must be one of the OT-FKM-type.
Therefore the cases for the multiplicity pairs $(4,5),(7,8)$ and
$(6,9)$ remain open now\footnote{Added in proof: Very recently Chi
\cite{Chi2} claimed a solution of the cases $(4,5),(6,9)$.}.

In the case of an isoparametric hypersurface with six principal
curvatures, M\"{u}nzner showed that all of the principal curvatures
must have the same multiplicity m, and Abresch \cite{A} showed that
$m$ must equal $1$ or $2$. By the classification of homogeneous
isoparametric hypersurfaces due to Takagi and Takahashi \cite{TT},
there is only one homogeneous family in each case up to congruence.
In the case of multiplicity $m = 1$, Dorfmeister and Neher \cite{DN}
showed that an isoparametric hypersurface must be homogeneous,
thereby completely classifying that case. It has long been
conjectured that the one homogeneous family in the case $g = 6$, $m
= 2$, is the only isoparametric family in this case, but this
conjecture has resisted proof for a long time. The approach that
Miyaoka \cite{Mi} used in the case $m = 1$ shows promise of possibly
leading to a proof of this conjecture, but so far a complete proof
has not been published\footnote{From private communications and
several geometric conferences, we learnt that R. Miyaoka
\cite{Miy09} had completed a proof of this conjecture.}.

While the classic theory of isoparametric hypersurfaces in real
space forms processed, people came to study similar objects in
complex space forms, \emph{i.e.}, hypersurfaces of constant
principal curvatures in $\mathbb{C}P^n$ or $\mathbb{C}H^n$ (cf.
\cite{Tak1}, \cite{Bern10} and references therein). In particular,
Berndt and D\'{\i}az-Ramos \cite{BD1, BD2} completely classified the
real hypersurfaces with $3$ distinct constant principal curvatures
in $\mathbb{C}H^n$ as that they are all homogeneous. Another
generalization is to consider global structures, \emph{i.e.}, a
family of parallel hypersurfaces of constant mean curvature, in
projective spaces $\mathbb{F}P^n$
($\mathbb{F}=\mathbb{C},\mathbb{H}$) just as \cite{Wa1} and
\cite{Pa} did, where \cite{Wa1} showed that there is a
correspondence between such ``isoparametric" hypersurfaces in
$\mathbb{C}P^n$ and isoparametric hypersurfaces in $S^{2n+1}$ by the
Hopf fibration, and there are some ``isoparametric" hypersurfaces in
$\mathbb{C}P^n$ with non-constant principal curvatures, while
\cite{Pa} investigated systematically the number of distinct
principal curvatures and their multiplicities of such
``isoparametric" hypersurfaces in $\mathbb{F}P^n$. Later, Terng and
Thorbergsson \cite{TeT} introduced a class of submanifolds in simply
connected symmetric spaces of compact type that they called
equifocal and similar structural results as the classic case were
established. A thorough study of the possible values of the
multiplicities $(m_1,m_2)$ for equifocal hypersurfaces was done by
Tang in \cite{Ta2}. It turns out that equifocal hypersurfaces are
the same as the level hypersurfaces of transnormal functions defined
in \cite{Wa2} as following.
\begin{defn}
A non-constant smooth function $f: N\rightarrow \mathbb{R}$ defined
on a Riemannian manifold $N$ is called \emph{transnormal} if there
is a smooth\footnote{All results would go through for $C^2$ category
instead of $C^{\infty}$.} function
$b:\mathbb{R}\rightarrow\mathbb{R}$ such that
\begin{equation}\label{iso1}
|\nabla f|^2=b(f).
\end{equation}
If moreover there is a continuous function
$a:\mathbb{R}\rightarrow\mathbb{R}$ such that
\begin{equation}\label{iso2}
\triangle f=a(f),
\end{equation}
then $f$ is called \emph{isoparametric} .
\end{defn}
As the roles of Cartan-M\"{u}nzner equations
(\ref{eq1.1},\ref{eq1.2}) in the classic theory, equation
(\ref{iso1}) implies that the regular hypersurfaces $M_t:=f^{-1}(t)$
(where $t$ is a regular value of $f$) are parallel and (\ref{iso2})
implies that these parallel hypersurfaces have constant mean
curvatures. These regular level hypersurfaces $M_t:=f^{-1}(t)$ of an
isoparametric function $f$ are called \emph{isoparametric
hypersurfaces}. The preimage of the maximum (resp. minimum) of an
isoparametric (or transnormal) function $f$ is called the
\emph{focal variety} of $f$, denoted by $M_{+}$ (resp. $M_{-}$).

In these terminologies Wang \cite{Wa2} proved a fundamental
structural result that the focal varieties $M_{\pm}$ of a
transnormal function are smooth submanifolds and each regular level
hypersurface $M_t$ is a tube over either of $M_{\pm}$, which is also
compatible with the classic theory of isoparametric hypersurfaces in
space forms. Based on this result, the authors \cite{GT3} improved
the fundamental theory of isoparametric functions on general
Riemannian manifolds and studied the existence problem of
isoparametric (transnormal) functions on certain Riemannian
manifolds, especially on exotic spheres. In particular, we proved
the non-existence of a properly transnormal function\footnote{Where
by \emph{proper} we mean that the focal varieties of the transnormal
function have codimension not less than $2$.} on any exotic
$4$-sphere (if exist). On the other hand, we obtained isoparametric
examples in some Brieskorn varieties and also in each Milnor sphere.
Moreover, we constructed explicitly a properly transnormal but not
isoparametric function on the Gromoll-Meyer sphere with two points
as the focal varieties, which also differs much from the classic
case due to a result claimed in \cite{Wa2} without proof that
regular level hypersurfaces of a (properly) transnormal function on
$S^{n+1}$ (or $\mathbb{R}^{n+1}$, but not for $\mathbb{H}^{n+1}$)
are isoparametric (cf. \cite{Mi2} for a proof).

By using Fermi coordinates, we \cite{GT4} completely proved the
minimality of the focal submanifolds of an isoparametric function on
a complete Riemannian manifold (claimed in \cite{Wa2} without
proof), and obtained the same properties of principal curvatures of
the focal submanifolds of an isoparametric function on a complete
Riemannian manifold as the classic case under the additional
assumption that each isoparametric hypersurface has constant
principal curvatures.

There're still many open problems\footnote{Added in proof: Very
recently Ge, Tang and Yan \cite{GTY} introduced the notion of $k$-th
isoparametric hypersurfaces which give a filtration of isoparametric
hypersurfaces by the number of constant higher order mean
curvatures. Hence many related questions arise such as existence
problem or classification problem of $k$-th isoparametric
hypersurfaces in certain Riemannian manifolds.} in this direction,
such as the existence problem of an isoparametric function on the
Gromoll-Meyer sphere with two points as its focal varieties, the
possibilities of the numbers $(g,m_1,m_2)$ for isoparametric
hypersurfaces in exotic spheres, the relationship with the normal
holonomy theory, \emph{etc}.

\section{Interactions and Applications}\label{int-app}
To conclude the story, in this section we would like to discuss some
interactions between Chern conjecture and isoparametric
hypersurfaces, and also some applications.

First, by the classification of possible values for $(g,m_1,m_2)$
and the simple structures of the second fundamental forms of them,
the minimal isoparametric hypersurfaces\footnote{There is a unique
minimal isoparametric hypersurface in each family of isoparametric
hypersurfaces in $S^{n+1}$, which also holds in a Riemannian
manifold of positive Ricci curvature (cf. \cite{GT3,GT4}).} and the
focal submanifolds in $S^{n+1}$ are the best examples satisfying the
Chern conjecture (resp. Conjecture \ref{lu conj}) that in this set
$\sigma$ (resp. $\sigma+\lambda_2$) takes only finite number of
values since it depends only on the values of $(g,m_1,m_2)$ and $n$.
Conversely, although the studies of Chern conjecture related to many
geometric or algebraic inequalities, such as the DDVV inequality and
its algebraic version, few applications to the theory of
isoparametric hypersurfaces could be found in literature. Further,
one could consider deducing some Simons-type inequalities in certain
Riemannian manifolds more general than in spheres such that
isoparametric hypersurfaces in these manifolds take similar roles as
those in spheres.

As the theory of isoparametric hypersurfaces processes extensively
and fruitfully, it seems that deeper study of its applications
worths people's attentions. For instance, by moving frame method,
Peng and Tang \cite{PTa} obtained the Brouwer degrees of gradient
maps (looked as maps from $S^{n+1}$ to itself) of isoparametric
functions (Cartan polynomials) on $S^{n+1}$ and applied them to
harmonic maps between spheres; Ma and Ohnita \cite{MO} studied
geometry of compact Lagrangian submanifolds in complex hyperquadrics
from the viewpoint of the theory of isoparametric hypersurfaces in
spheres, where they determined the Hamiltonian stability of compact
minimal Lagrangian submanifolds embedded in complex hyperquadrics
which are obtained as Gauss images of isoparametric hypersurfaces in
spheres with $g(= 1, 2, 3)$ distinct principal curvatures; Ge and
Xie \cite{GX10} studied properties of the gradient map of an
isoparametric function (Cartan polynomial) which help them to deduce
the Brouwer degrees of gradient maps immediately and to construct a
counterexample to the Br\'ezis question (cf. \cite{Br}) on the
symmetry for the Ginzburg-Landau system in dimension $6$, which
gives a partial answer towards the Open problem $2$ raised by Farina
\cite{Far}. Note that the theory had also affections in differential
topology as we observed in \cite{GT3}, and it is obviously of great
interest that further applications can be found.


\end{document}